\documentclass[11pt]{amsart}

\newtheorem{proposition}{Proposition}[section]
\newtheorem{lemma}[proposition]{Lemma}

\newtheorem{theorem}[proposition]{Theorem}

\theoremstyle{definition}
\newtheorem{definition}[proposition]{Definition}

\theoremstyle{remark}

\newcommand{\thlabel}[1]{\label{th:#1}}
\newcommand{\thref}[1]{Theorem~\ref{th:#1}}
\newcommand{\selabel}[1]{\label{se:#1}}
\newcommand{\seref}[1]{Section~\ref{se:#1}}
\newcommand{\lelabel}[1]{\label{le:#1}}

\newcommand{\prlabel}[1]{\label{pr:#1}}
\newcommand{\prref}[1]{Proposition~\ref{pr:#1}}

\newcommand{\delabel}[1]{\label{de:#1}}

\newcommand{\eqlabel}[1]{\label{eq:#1}}
\newcommand{\equref}[1]{(\ref{eq:#1})}

\newcommand{\Hom}{{\rm Hom}}

\newcommand{\End}{{\rm End}}

\newcommand{\can}{{\rm can}}

\def\ot{\otimes}

\newcommand{\Cc}{\mathcal{C}}
\newcommand{\Dd}{\mathcal{D}}

\newcommand{\Mm}{\mathcal{M}}

\def\text#1{{\rm {\rm #1}}}

\def\ol{\overline}

\begin{document}
\title[Galois corings applied to partial Galois theory]{Galois corings applied to partial Galois theory}
\author{S. Caenepeel}
\address{Faculty of Engineering Sciences,
Vrije Universiteit Brussel, VUB, B-1050 Brussels, Belgium}
\email{scaenepe@vub.ac.be}
\urladdr{http://homepages.vub.ac.be/\~{}scaenepe/}
\author{E. De Groot}
\address{Faculty of Engineering Sciences,
Vrije Universiteit Brussel, VUB, B-1050 Brussels, Belgium}
\email{edegroot@vub.ac.be}
\urladdr{http://homepages.vub.ac.be/\~{}edegroot/}
\begin{abstract}
Partial Galois extensions were recently introduced by Doku\-chaev, Ferrero
and Paques. We introduce partial Galois extensions for noncommutative rings,
using the theory of Galois corings. We associate a Morita context to a
partial action on a ring.
\end{abstract}

\maketitle

\section*{Introduction}
Partial actions of groups originate from the theory of operator algebras, see for example \cite{Exel}. Partial representations of groups on Hilbert spaces were introduced
independently in \cite{Exel2} and \cite{Raeburn}. Several applications are given
in the literature, we refer to \cite{DFP} for a more extensive bibliography.
More recently, partial actions were studied from a purely algebraic point of
view, in \cite{DEP,DE,DZ}.\\
In \cite{DFP}, the authors consider partial actions on commutative rings,
with the additional assumption that the associated ideals are generated by
idempotents. Then they generalize Galois theory for commutative rings,
as introduced in \cite{CHR} for usual group actions, to partial actions.\\
Corings were introduced by Sweedler in 1975 in \cite{S3}. There has been
a revived interest in corings since the beginning of the century, based on
an observation made by Takeuchi that various types of modules, such as
Hopf modules, relative Hopf modules, graded modules, entwined modules and
Yetter-Drinfeld modules may be viewed as comodules over a coring. 
Brzezi\'nski \cite{B1} noticed the importance of this observation: the language
of corings can be applied successfully to give a unified and more elegant
treatment to properties related to all these kinds of modules. An overview
can be found in \cite{BW}.\\
One of the nice applications is descent and Galois theory: Galois corings
were introduced in \cite{B1}, and studied in \cite{Caenepeel03} and
\cite{W3}. The corings approach provides a unified theory for various
types of Galois theories, including the classical Chase-Harrison-Rosenberg theory
\cite{CHR}, Hopf-Galois theory (see \cite{CS,KT,Schneider90}),
coalgebra Galois theory (see \cite{BrzezinskiH}) and weak Hopf-Galois
theory (see the forthcoming \cite{CDG2}).\\
The aim of this note is to develop partial Galois theory starting from
Galois corings. The strategy is basically the following: given a set of
idempotents $e_\sigma$ indexed by a finite group $G$ in a ring $A$, we investigate when the direct sum of the $Ae_\sigma$ is a coring; it turns out that this is the
case if a partial action of $G$ on $A$ is given. Then we investigate when
this coring is a Galois coring, and apply the results in \cite{Caenepeel03}.
This procedure still works in the case where the ring $A$ is not commutative.
In the case where $A$ is commutative, we
recover some of the results in \cite{DFP}. This is done in \seref{4}.
In \seref{3}, we associate a Morita context to a partial action on a ring
$A$, and show that the context is strict if $A$ is a faithfully flat partial
Galois extension of the invariants ring $A^G$.   

\section{Preliminary results}\selabel{1}
\subsection{Galois corings}\selabel{1.1}
Let $A$ be a ring. An $A$-coring $\Cc$ is a coalgebra in the category
${}_A\Mm_A$ of $A$-bimodules. Thus an $A$-coring
is a triple $\Cc=(\Cc,\Delta_\Cc,\varepsilon_\Cc)$, where
 $\Cc$ is an $A$-bimodule, and
$\Delta_\Cc:\ \Cc\to \Cc\ot_A\Cc$  and
 $\varepsilon_\Cc:\ \Cc\to A$ are $A$-bimodule maps
such that
\begin{equation}\eqlabel{1.1.1}
(\Delta_\Cc\ot_A \Cc)\circ \Delta_\Cc=
(\Cc\ot_A\Delta_\Cc )\circ \Delta_\Cc,
\end{equation}
and
\begin{equation}\eqlabel{1.1.2}
({\Cc}\ot_A\varepsilon_{\Cc})\circ \Delta_\Cc=
(\varepsilon_{\Cc}\ot_A{\Cc})\circ \Delta_\Cc=\Cc.
\end{equation}
We use the Sweedler-Heyneman notation for the comultiplication:
$$\Delta_\Cc(c)=c_{(1)}\ot_A c_{(2)}.$$
A right $\Cc$-comodule $M=(M,\rho )$ consists of a right $A$-module
$M$ together with a right $A$-linear map $\rho :\ M\to M\ot_A\Cc$
such that:
\begin{equation}\eqlabel{1.1.3}
(\rho \ot_A \Cc)\circ \rho =(M\ot_A\Delta_\Cc)\circ \rho, 
\end{equation}
and
\begin{equation}\eqlabel{1.1.4}
(M\ot_A\varepsilon_\Cc)\circ \rho =M.
\end{equation}
We then say that $\Cc$ coacts from the right on $M$, and we denote
$$\rho (m)=m_{[0]}\ot_A m_{[1]}.$$
A right $A$-linear map
$f:\ M\to N$ between two right $\Cc$-comodules $M$ and $N$ is called
right $\Cc$-colinear if $\rho (f(m))=f(m_{[0]})\ot m_{[1]}$,
for all $m\in M$. The category of right $\Cc$-comodules and $\Cc$-colinear
maps is denoted by $\Mm^\Cc$.\\
$x\in \Cc$ is called grouplike if $\Delta_\Cc(x)=x\ot x$ and
$\varepsilon_\Cc(x)=1$. Grouplike elements of $\Cc$ correspond bijectively
to right $\Cc$-coactions on $A$: if $A$ is grouplike, then we have the following
right $\Cc$-coaction $\rho$ on $A$: $\rho(a)=xa$.\\
Let $(\Cc,x)$ be a coring with a fixed grouplike element. For $M\in \Mm^{\Cc}$, we call
$$M^{{\rm co}\Cc}=\{m\in M~|~\rho (m)=m\ot_A x\}$$
the submodule of coinvariants of $M$. Observe that
$$A^{{\rm co}\Cc}=\{b\in A~|~bx=xb\}$$
is a subring of $A$. Let $i:\ B\to A$ be a ring morphism.
$i$ factorizes through $A^{{\rm co}\Cc}$ if and only if
$$x\in G(\Cc)^B=\{x\in G(\Cc)~|~xb=bx,~{\rm for~all~}b\in B\}.$$ 
We then have
a pair of adjoint functors
$(F,G)$, respectively between the categories $\Mm_B$ and
$\Mm^\Cc$ and the categories ${}_B\Mm$ and
${}^\Cc\Mm$. For 
$N\in \Mm_B$ and $M\in \Mm^\Cc$,
$$F(N)=N\ot_B A~~{\rm and}~~G(M)=M^{{\rm co}\Cc}.$$
The unit and counit of the adjunction are
$$\nu_N:\ N\to (N\ot_BA)^{{\rm co}\Cc},~~\nu_N(n)=n\ot_B1;$$
$$\zeta_M:\ M^{{\rm co}\Cc}\ot_B A\to M,~~
\zeta_M(m\ot_B a)=ma.$$
Let $i:\ B\to A$ be a morphism of rings. The associated canonical coring is
$\Dd=A\ot_BA$, with comultiplication and counit given by the formulas
$$\Delta_\Dd:\ \Dd\to \Dd\ot_A\Dd\cong A\ot_B A\ot_B A,~~
\Delta_\Dd(a\ot_B a')=a\ot_B1\ot_Ba'$$
and
$$\varepsilon_\Dd:\ \Dd=A\ot_B A\to A,~~\varepsilon_\Dd(a\ot_Ba')=aa'.$$
If $i:\ B\to A$ is pure as a morphism of left and right $B$-modules, then
the categories $\Mm_B$ and $\Mm^\Dd$ are equivalent.\\
Let $(\Cc,x)$ be a coring with a fixed grouplike element, and $i:\ B\to A^{{\rm co}\Cc}$
a ring morphism. We then have a morphism of corings
$${\rm can}:\ \Dd=A\ot_BA\to \Cc,~~{\rm can}(a\ot_B a')=axa'.$$
If $F$ is fully faithful, then $B\cong A^{{\rm co}\Cc}$; if $G$ is fully
faithful, then ${\rm can}$ is  an isomorphism.
$(\Cc,x)$ is called a Galois coring if  
${\rm can}:\ A\to_{A^{{\rm co}\Cc}}A\to \Cc$ is
bijective. From \cite{Caenepeel03}, we recall the following results.

\begin{theorem}\thlabel{1.2}
Let $(\Cc,x)$ be an $A$-coring with fixed grouplike element, and
$B=A^{{\rm co}\Cc}$. Then the following statements are equivalent.
\begin{enumerate}
\item $(\Cc,x)$ is Galois and $A$ is faithfully flat as a left $B$-module;
\item $(F,G)$ is an equivalence
and $A$ is flat as a left $B$-module.
\end{enumerate}
\end{theorem}

Let $(\Cc,x)$ be a coring with a fixed grouplike element, and take
$T=A^{{\rm co}\Cc}$. Then ${}^*\Cc={}_A\Hom(\Cc,A)$ is a ring, with
multiplication given by
\begin{equation}\eqlabel{1.2.1a}
(f\# g)(c)=g(c_{(1)}f(c_{(2)})).
\end{equation}
We have a morphism of rings $j:\ A\to {}^*\Cc$, given by
$$j(a)(c)=\varepsilon_\Cc(c)a.$$
This makes ${}^*\Cc$ into an $A$-bimodule,
via the formula
$$(afb)(c)=f(ca)b.$$
Consider the left dual of the canonical map:
$${}^*{\rm can}:\ {}^*\Cc\to {}^*\Dd\cong {}_T\End(A)^{\rm op},~~
{}^*{\rm can}(f)(a)=f(xa).$$
We then have the following result.

\begin{proposition}\prlabel{2.2}
If $(\Cc,x)$ is Galois, then ${}^*{\rm can}$ is an isomorphism.
The converse property holds if $\Cc$ and $A$ are finitely generated
projective, respectively as a left $A$-module, and a left $T$-module.
\end{proposition}

Let $Q=\{q\in {}^*\Cc~|~c_{(1)}q(c_{(2)})=q(c)x,~{\rm for~all~}c\in \Cc\}$.
A straightforward computation shows that $Q$ is a $({}^*\Cc,T)$-bimodule.
Also $A$ is a left $(T,{}^*\Cc)$-bimodule; the right ${}^*\Cc$-action is
induced by the right $\Cc$-coaction: $a\cdot f=f(xa)$. Now
consider the maps
\begin{eqnarray}
\tau:&& A\ot_{{}^*\Cc} Q\to T,~~\tau(a\ot_{{}^*\Cc} q)=q(xa);\eqlabel{2.2.1}\\
\mu:&& Q\ot_T A\to {}^*\Cc,~~\mu(q\ot_T a)=q\#i(a).\eqlabel{2.2.2}
\end{eqnarray}
With this notation, we have the following property (see \cite{CVW}).

\begin{proposition}\prlabel{2.3}
$(T,{}^*\Cc,A,Q,\tau,\mu)$ is a Morita context.
\end{proposition}

We also have (see \cite{Caenepeel03}):

\begin{theorem}\thlabel{2.4}
Let $(\Cc,x)$ be a coring with fixed grouplike element, and assume
that $\Cc$ is a left $A$-progenerator. We take a subring
$B$ of $T=A^{{\rm co}\Cc}$, and consider the map
$${\rm can}:\ \Dd=A\ot_{B}A\to \Cc,~~{\rm can}(a\ot_{T}a')=axa'$$
Then the following statements are equivalent:
\begin{enumerate}
\item \begin{itemize}
\item ${\rm can}$ is an isomorphism;
\item $A$ is faithfully flat as a left $B$-module.
\end{itemize}
\item \begin{itemize}
\item ${}^*{\rm can}$ is an isomorphism;
\item $A$ is a left $B$-progenerator.
\end{itemize}
\item \begin{itemize}
\item $B=T$;
\item the Morita context $(B,{}^*\Cc,A,Q,\tau,\mu)$ is strict.
\end{itemize}
\item \begin{itemize}
\item $B=T$;
\item $(F,G)$ is an equivalence of categories.
\end{itemize}
\end{enumerate}
\end{theorem}

\subsection{Partial group actions}\selabel{1.3}
Let $G$ be a finite group, and $R\to S$ a commutative ring extension.
From \cite{DE}, we recall that a partial action $\alpha$ of $G$ on $S$
is a collection of ideals $S_\sigma$ and isomorphisms of ideals
$\alpha_\sigma:\ S_{\sigma^{-1}}\to S_\sigma$ such that
\begin{enumerate}
\item $S_1=S$, and $\alpha_1=S$, the identity on $S$;
\item $S_{(\sigma\tau)^{-1}}\supset\alpha_\tau^{-1}(S_\tau\cap S_{\sigma^{-1}})$;
\item $(\alpha_\sigma\circ \alpha_\tau)(x)=\alpha_{\sigma\tau}(x)$, for all
$x\in \alpha_\tau^{-1}(S_\tau\cap S_{\sigma^{-1}})$.
\end{enumerate}
In \cite{DFP}, the following particular situation is considered: every
$S_\sigma$ is of the form $S_\sigma=Se_{\sigma}$, where $e_{\sigma}$ is an
idempotent of $S$. In this case, we can show that
\begin{equation}\eqlabel{1.3.1}
\alpha_\sigma(\alpha_\tau(xe_{\tau^{-1}})e_{\sigma^{-1}})=\alpha_{\sigma\tau}(
xe_{\tau^{-1}\sigma^{-1}})e_\sigma,
\end{equation}
for all $\sigma,\tau\in G$ and $x\in S$. We then have an associative ring
with unit 
$$A\star_\alpha G=\bigoplus_{\sigma\in G} Ae_\sigma u\sigma,$$
with multiplication
\begin{equation}\eqlabel{1.3.2}
(a_\sigma u_\sigma)(a_\tau u_\tau)=\alpha_\sigma(\alpha_{\sigma^{-1}}(a_\sigma)
b_\tau)u_{\sigma\tau}.
\end{equation}

\section{Partial Galois theory for noncommutative rings}\selabel{4}
Let $A$ be a (noncommutative) ring, and $G$ a finite group. For every
$\sigma\in G$, we assume that there is a central idempotent $e_\sigma\in A$,
and a ring automorphism
$$\alpha_\sigma:\ Ae_{\sigma^{-1}}\to Ae_\sigma.$$
In particular, it follows that $\alpha_\sigma(e_{\sigma^{-1}})=e_{\sigma}$.
We can extend $\alpha_\sigma$ to $A$, by putting $\alpha_\sigma(a)=
\alpha_\sigma(ae_\sigma)$, for all $a\in A$.\\
Then we consider the direct sum $\Cc$ of all the $Ae_{\sigma}$. Let $v_\sigma$
be the element of $\Cc$ with $e_\sigma$ in the $Ae_\sigma$-component, and
$0$ elsewhere. We then have
$$\Cc=\bigoplus_{\sigma\in G} Ae_{\sigma}v_{\sigma}=\bigoplus_{\sigma\in G} Av_{\sigma}.$$
Obviously $\Cc$ is a left $A$-module.

\begin{lemma}\lelabel{4.1}
$\Cc$ is an $A$-bimodule. The right $A$-action is given by the formula
\begin{equation}\eqlabel{4.1.1}
(a'v_\sigma)a=a'\alpha_\sigma(ae_{\sigma^{-1}})v_\sigma
\end{equation}
\end{lemma}

\begin{proof}
Let us show that \equref{4.1.1} is an associative action: for all $a,a'\in A$,
we have
\begin{eqnarray*}
&&\hspace*{-2cm}
v_\sigma(aa')=\alpha_\sigma(aa'e_{\sigma^{-1}})v_\sigma=
\alpha_\sigma(ae_{\sigma^{-1}}a'e_{\sigma^{-1}})v_\sigma\\
&=& \alpha_\sigma(ae_{\sigma^{-1}})\alpha_\sigma(a'e_{\sigma^{-1}})v_\sigma=
\alpha_\sigma(ae_{\sigma^{-1}})v_\sigma a'= (v_\sigma a)a'.
\end{eqnarray*}
\end{proof}

We now consider the left $A$-linear maps
$$\Delta_\Cc:\ \Cc\to \Cc\ot_A \Cc,~~\Delta_\Cc(av_\sigma)=
\sum_{\tau\in G} av_\tau\ot_A v_{\tau^{-1}\sigma};$$
$$\varepsilon_\Cc:\ \Cc\to A,~~\varepsilon_\Cc(\sum_{\sigma\in G}a_\sigma v_\sigma)=
a_1.$$

\begin{proposition}\prlabel{4.2}
With notation as above, $(\Cc,\Delta_\Cc,\varepsilon_\Cc)$ is an $A$-coring if
and only if $e_1=A$, $\alpha_1=A$ and
\begin{equation}\eqlabel{4.2.1}
\alpha_\sigma(\alpha_\tau(ae_{\tau^{-1}})e_{\sigma^{-1}})=
\alpha_{\sigma\tau}(ae_{\tau^{-1}\sigma^{-1}})e_\sigma,
\end{equation}
for all $a\in A$ and $\sigma,\tau\in G$.
\end{proposition}

\begin{proof}
We compute
\begin{eqnarray*}
\Delta_\Cc(v_\sigma a)&=&
\Delta_\Cc(\alpha_\sigma(ae_{\sigma^{-1}})v_\sigma\\
&=& \sum_{\tau\in G} \alpha_\sigma(ae_{\sigma^{-1}}) v_\tau\ot_A v_{\tau^{-1}\sigma};\\
\Delta_\Cc(v_\sigma)a&=&
\sum_{\tau\in G}v_\tau\ot_A v_{\tau^{-1}\sigma} a\\
&=& \sum_{\tau\in G} v_\tau \ot_A \alpha_{\tau^{-1}\sigma}(ae_{\sigma^{-1}\tau})
v_{\tau^{-1}\sigma}\\
&=& \sum_{\tau\in G} \alpha_\tau(\alpha_{\tau^{-1}\sigma}(ae_{\sigma^{-1}\tau})
e_{\tau^{-1}})v_\tau\ot_A v_{\tau^{-1}\sigma}.
\end{eqnarray*}
Hence $\Delta_\Cc$ is right $A$-linear if and only if
$$\alpha_\sigma(ae_{\sigma^{-1}})e_\tau =
\alpha_\tau(\alpha_{\tau^{-1}\sigma}(ae_{\sigma^{-1}\tau})
e_{\tau^{-1}}),$$
for all $\sigma,\tau\in G$ and $a\in A$. Substituting $\lambda=\tau^{-1}\sigma$,
we find that this is equivalent to \equref{4.2.1}.\\
Let us now investigate when $\varepsilon_\Cc$ is right $A$-linear. We have
$$\varepsilon_\Cc(\sum_{\sigma\in G} a_\sigma v_\sigma)a=a_1a)$$
and
$$\varepsilon_\Cc(\sum_{\sigma\in G} a_\sigma v_\sigma a)=
\varepsilon_\Cc(\sum_{\sigma\in G} a_\sigma \alpha_\sigma(ae_{\sigma^{-1}})v_\sigma)
=a_1\alpha_1(a e_1)$$
If $\varepsilon_\Cc$ is right $A$-linear, then we find that
$\alpha_1(a e_1)=a$, for all $a\in A$. In particular,
$$e_1= \alpha_1(e_1 e_1)= \alpha_1( 1 e_1)=1,$$
and then it follows that $\alpha_1(a)=a$, for all $a\in A$.
Conversely, if $e_1=1$ and $\alpha_1=A$, then it follows that $\varepsilon_\Cc$
is right $A$-linear.\\
Now assume that \equref{4.2.1} holds, and that $e_1=1$ and $\alpha_1=A$.
The coassociativity and counit property then follow in a straightforward way.
\end{proof}

From now on, we will assume that $\Cc=\bigoplus_{\sigma\in G} Av_\sigma$ is
an $A$-coring. The set of data $(e_\sigma,\alpha_\sigma)_{\sigma\in G}$ will
be called an idempotent partial action of $G$ on $A$. This is the case for
the partial actions discussed in \cite{DFP}, that we recalled in \seref{1.3},
in view of \equref{1.3.1}.

\begin{lemma}\lelabel{4.3}
$x=\sum_{\sigma\in G} v_\sigma$ is  a grouplike element of $\Cc$.
\end{lemma}

\begin{proof}
$\varepsilon_\Cc(x)=1$, and
$$
\Delta_\Cc(x)=\sum_{\sigma,\tau\in G} v_\tau\ot_Av_{\tau^{-1}\sigma}=
\sum_{\rho,\tau\in G} v_\tau\ot_Av_{\rho}=x\ot_A x.$$
\end{proof}

Consider the left $A$-linear maps
\begin{equation}\eqlabel{4.4.0}
u_\sigma:\ \Cc\to Ae_\sigma,~~u_\sigma(\sum_{\tau\in G} a_\tau v_\tau)=a_\sigma e_\sigma.
\end{equation}
Then for all $c\in \Cc$, we have
$$c=\sum_{\sigma\in G} u_{\sigma}(c)v_{\sigma},$$
hence $\{(u_\sigma,v_\sigma)~|~\sigma\in G\}$ is a dual basis of $\Cc$ as a
left $A$-module.\\
Now let $(M,\rho)$ be a right $\Cc$-comodule. We have a right $A$-linear map
$$\rho:\ M\to \bigoplus_{\sigma\in G} M\ot_A Av_{\sigma}.$$
Consider the maps
$$\rho_\sigma=(M\ot_A u_\sigma)\circ \rho:\ M\to Me_\sigma.$$
We then have
$$\rho(m)=m_{[0]}\ot_A m_{[1]}=\sum_{\sigma\in G}\rho_\sigma(m)\ot_A v_\sigma,$$
for all $m\in M$. From the fact $\rho$ is right $A$-linear, it follows that
$$\rho(ma)=\sum_{\sigma\in G} \rho_\sigma(ma)\ot_A v_\sigma=
\rho(m)a= \sum_{\sigma\in G} \rho_\sigma(m)\ot_A \alpha_\sigma(ae_{\sigma^{-1}})v_\sigma,
$$
hence
\begin{equation}\eqlabel{4.4.1}
\rho_{\sigma}(ma)=\rho_{\sigma}(m)\alpha_\sigma(ae_{\sigma^{-1}}),
\end{equation}
for all $m\in M$ and $\sigma\in G$. It follows from \equref{4.4.1} that
$$\rho_\sigma(me_{\sigma^{-1}})=\rho_\sigma(m)\alpha_\sigma(e_{\sigma^{-1}})=
\rho_\sigma(m)e_\sigma.$$
This means that $\rho_\sigma:\ M\to Me_\sigma$ factors through the projection
$M\to M e_{\sigma^{-1}}$, so we obtain a map
$$\rho_\sigma:\ Me_{\sigma^{-1}}\to Me_\sigma.$$
Since $(M\ot_A\varepsilon_\Cc)\circ \rho=M$, we have, for all $m\in M$:
$$m=\sum_{\sigma\in G}\rho_\sigma(m)\ot_A \varepsilon_\Cc(v_\sigma)=
\rho_1(m) e_1=\rho_1(m).$$
Hence $\rho_1:\ Me_1=M\to Me_1=M$ is the identity. From the coassociativity
of $\rho$, we deduce that
\begin{eqnarray*}
&&\hspace*{-2cm}
\sum_{\sigma,\tau\in G}\rho_\tau(\rho_\sigma(m))\ot_A v_\tau\ot_A v_\sigma=
\sum_{\sigma,\rho\in G} \rho_\sigma(m)\ot_A v_\mu\ot_A v_{\mu^{-1}\sigma}\\
&=& \sum_{\kappa,\mu\in G}v_{\mu\circ\kappa}\ot_A v_\mu\ot_A v_\kappa
=\sum_{\sigma,\tau\in G}v_{\tau\circ\sigma}\ot_A v_\tau\ot_A v_\sigma,
\end{eqnarray*}
hence
\begin{equation}\eqlabel{4.4.2}
\rho_\tau(\rho_\sigma(m))=\rho_{\tau\circ\sigma}(m),
\end{equation}
for all $m\in M$ and $\sigma,\tau\in G$. In particular,
\begin{equation}\eqlabel{4.4.3}
\rho_\tau(\rho_\sigma(me_{\sigma^{-1}}e_{\tau^{-1}}))=
\rho_{\tau\circ \sigma}(me_{\sigma^{-1}\tau^{-1}})e_\tau.
\end{equation}
It follows from \equref{4.4.3} that
$\rho_{\tau^{-1}}:\ Me_\tau\to Me_{\tau^{-1}}$ is the inverse of
$\rho_{\tau}:\ Me_{\tau^{-1}}\to Me_\tau$.

\begin{definition}\delabel{4.4}
Let $(e_\sigma,\alpha_\sigma)_{\sigma\in G}$ be an idempotent partial action
of $G$ on $A$, and $M$ a right $A$-module. A partial Galois descent datum
consists of a set of maps
$$\rho_\sigma:\ M\to Me_\tau$$
such that $\rho_1=M$, the identity on $M$, the restriction of
$\rho_\sigma$ to $Me_{\tau^{-1}}$ is an isomorphism, and
\equref{4.4.1} and \equref{4.4.3} hold for all $m\in M$, $a\in A$
and $\sigma,\tau\in G$.
\end{definition}

\begin{proposition}\prlabel{4.5}
Let $(e_\sigma,\alpha_\sigma)_{\sigma\in G}$ be an idempotent partial action
of $G$ on $A$, and $\Cc$ the corresponding $A$-coring. Then right $\Cc$-coactions
on $M$ correspond bijectively to partial Galois descent data.
\end{proposition}

\begin{proof}
We have already explained above how a right $\Cc$-coaction $\rho$ on $M$ can be
transformed into a partial Galois descent datum. Conversely, let
$(\rho_\sigma)_{\sigma\in G}$ be a partial Galois descent datum, and
define $\rho:\ M\to M\ot_A\Cc$ by
$$\rho(m)=\sum_{\sigma\in G} \rho_{\sigma}(m)\ot_A v_\sigma.$$
Straightforward computations show that $\rho$ is a coaction, and that the two
constructions are inverse to each other.
\end{proof}

Let $M$ be a right $\Cc$-comodule. Then $m\in M^{{\rm co}\Cc}$ if and only if
$$\rho(m)=\sum_{\sigma\in G} \rho_\sigma(m)\ot_A v_\sigma=
\sum_{\sigma\in G} m\ot_A v_\sigma= \sum_{\sigma\in G} me_\sigma\ot_A v_\sigma$$
if and only if
$$\rho_\sigma(m)=\rho_\sigma(me_{\sigma^{-1}})=me_\sigma,$$
for all $\sigma\in G$. We define
$$M^G=\{m\in M~|~\rho_\sigma(me_{\sigma^{-1}})=me_{\sigma^{-1}},~
{\rm for~all~}\sigma\in G\}=M^{\Cc}.$$
The grouplike element $x=\sum_{\sigma\in G}v_\sigma$ makes $A$ into a right
$\Cc$-comodule:
$$\rho(a)=1\ot_A xa=\sum_{\sigma\in G}\alpha_\sigma(ae_{\sigma^{-1}})\ot_Av_\sigma,$$
and we have
$$T=A^G=\{a\in A~|~\alpha_\sigma(ae_{\sigma^{-1}})=ae_\sigma,~
{\rm for~all~}\sigma\in G\}.$$
Let $i:\ B\to T$ be a ring morphism.
We have seen in \seref{1.1} that we have a pair of adjoint functors $(F,G)$:
$$F:\ \Mm_B\to \Mm^\Cc,~~F(N)=N\ot_B A;$$
$$G:\ \Mm^\Cc\to \Mm_B,~~G(N)=N^G.$$
$F(N)=N\ot_B A$ is a right $\Cc$-comodule in the following way:
$$\rho_\sigma(n\ot_A a)=n\ot \alpha_\sigma(a).$$
The canonical map is the following:
$${\rm can}:\ A\ot_B A\to \bigoplus_{\sigma\in G}A_ev_\sigma,~~
{\rm can}(a\ot b)=\sum_{\sigma\in G}a\alpha_\sigma(be_{\sigma^{-1}})v_\sigma.$$
$\bigoplus_{\sigma\in G}A_ev_\sigma$ is a Galois coring if ${\rm can}:
A\ot_{A^G} A\to \bigoplus_{\sigma\in G}A_ev_\sigma$ is
an isomorphism. We will then say that $A$ is a partial $G$-Galois extension
of $A^G$. From \thref{1.2}, we immediately obtain the following result.

\begin{theorem}\thlabel{4.5a}
Let $(e_\sigma,\alpha_\sigma)_{\sigma\in G}$ be an idempotent partial action
of $G$ on $A$, and $T=A^G$. Then the following assertions are equivalent.
\begin{enumerate}
\item  $A$ is a partial $G$-Galois extension
of $T$ and $T$ is faithfully flat as a left $T$-module;
\item $(F,G)$ is a category equivalence and $A$ is flat as a left $T$-module.
\end{enumerate}
\end{theorem}

\section{Partial actions and Morita theory}\selabel{3}
Let us now compute the multiplication on 
$${}^*\Cc={}_A\Hom(\Cc,A)=\bigoplus_{\sigma\in G} {}_A\Hom(Ae_\sigma,A).$$
 We will use
the maps $u_\sigma$ defined in \equref{4.4.0}. Also recall that ${}^*\Cc$
is an $A$-bimodule, with left and right $A$-action
$$(afb)(c)=f(ca)b.$$
Take $f\in {}_A\Hom(\Cc, A)$. For all $c\in\Cc$, we have
$$f(c)=\sum_{\sigma\in \Cc} u_\sigma(c)f(v_\sigma).$$
Now $f(v_\sigma)=f(e_\sigma v_\sigma)=e_\sigma f(v_\sigma)\in Ae_\sigma$,
so we can conclude that
$${}^*\Cc=\bigoplus_{\sigma\in  \Cc} u_\sigma Ae_\sigma.$$
For $b\in A$, we compute
$$
(bu_\tau)(\sum_{\sigma\in G} a_\sigma v_\sigma)=
u_\tau(\sum_{\sigma\in G} a_\sigma v_\sigma b)=
u_\tau(\sum_{\sigma\in G} a_\sigma\alpha_\sigma(be_{\sigma^{-1}})v_\sigma)=
a_\tau\alpha_\tau(be_{\tau^{-1}}),$$
and we conclude that
\begin{equation}\eqlabel{4.6.1}
bu_\tau=u_\tau\alpha_\tau(be_{\tau^{-1}}).
\end{equation}
We next compute, using \equref{4.2.1}: 
\begin{eqnarray*}
&&\hspace*{-2cm}(u_\rho\# u_\nu)(\sum_{\sigma\in G} a_\sigma v_\sigma)= u_\nu(\sum_{\sigma,\tau}a_\sigma v_\tau u_\rho(v_{\tau^{-1}\sigma}))
= u_\nu(\sum_{\sigma,\tau}a_\sigma v_\tau \delta_{\tau\rho,\sigma})\\
&=&
 u_\nu(\sum_\tau a_{\tau\rho}v_\tau)=a_{\nu\rho}=
  u_{\nu\rho}(\sum_{\sigma\in G} a_\sigma v_\sigma),
\end{eqnarray*}
and we conclude that
\begin{equation}\eqlabel{4.6.2}
u_\sigma\# u_\tau =u_{\sigma\tau}.
\end{equation}
We can summarize this as follows:

\begin{proposition}\prlabel{4.6}
Let $\Cc=\bigoplus_{\sigma\in G}Av_\sigma$. The left dual ring is
$${}^*\Cc=\bigoplus_{\sigma\in G}u_\sigma Ae_\sigma,$$
with multiplication rule
\begin{equation}\eqlabel{4.6.3}
u_\tau b_\tau\# u_\sigma a_\sigma= 
u_{\sigma\tau}\alpha_\sigma(b_\tau e_{\sigma^{-1}})a_{\sigma}.
\end{equation}
\end{proposition}

If $A$ is commutative, then ${}^*\Cc$ is isomorphic to $(A\star_\alpha G)^{\rm op}$,
as introduced in \cite{DFP}, see \equref{1.3.2}. Indeed, for
$a\in Ae_\sigma$ and $b\in Ae_\tau$, we compute that
\begin{eqnarray*}
&&\hspace*{-2cm}
\alpha_\sigma(\alpha_{\sigma^{-1}}(a_\sigma)b_\tau)=
\alpha_\sigma(\alpha_{\sigma^{-1}}(a_\sigma)b_\tau e_{\sigma^{-1}})\\
&=& a_\sigma \alpha_\sigma(b_\tau e_{\sigma^{-1}})
\alpha_\sigma(b_\tau e_{\sigma^{-1}})a_{\sigma}
\end{eqnarray*}

Recall that a ring morphism $A\to R$ is called {\sl Frobenius}
if there exists an $A$-bimodule map $\ol{\nu}:\ R\to A$ and
$e=e^1\ot_A e^2\in R\ot_AR$ (summation implicitly understood) such that
\begin{equation}\eqlabel{4.7.1}
re^1\ot_A e^2=e^1\ot_A e^2r
\end{equation}
for all $r\in R$, and
\begin{equation}\eqlabel{4.7.2}
\ol{\nu}(e^1)e^2=e^1\ol{\nu}(e^2)=1.
\end{equation}
This is equivalent to the restrictions of scalars $\Mm_R\to
\Mm_A$ being a Frobenius functor, which means that its left and right
adjoints are isomorphic (see \cite[Sec. 3.1 and 3.2]{CMZ}).
$(e,\ol{\nu})$ is then called a Frobenius system.

\begin{proposition}\prlabel{4.7}
Suppose that we have an idempotent partial action of $G$ on $A$. Then
the ring morphism $A\to {}^*\Cc$ is Frobenius.
\end{proposition}

\begin{proof}
The Frobenius system is
$(e=\sum_{\sigma\in G} u_{\sigma^{-1}}\ot_A u_\sigma, \ol{\nu})$,
with
$$\ol{\nu}(\sum_{\sigma\in G} u_\sigma a_{\sigma})=a_1.$$
We compute that, for all $a\in A$,
\begin{eqnarray*}
&&\hspace*{-2cm}
a\sum_{\sigma\in G} u_{\sigma^{-1}}\ot_A u_\sigma=
\sum_{\sigma\in G} u_{\sigma^{-1}}\alpha_{\sigma^{-1}}(ae_\sigma)\ot_A u_\sigma\\
&=& \sum_{\sigma\in G} u_{\sigma^{-1}}\ot_A u_\sigma 
\alpha_\sigma(\alpha_{\sigma^{-1}}(ae_\sigma))=
\sum_{\sigma\in G} u_{\sigma^{-1}}\ot_A u_\sigma ae_\sigma\\
&=& \sum_{\sigma\in G} u_{\sigma^{-1}}\ot_A u_\sigma a.
\end{eqnarray*}
The rest is obvious.
\end{proof}

Let $i:\ B\to T=A^{{\rm co}\Cc}$ be a ring morphism. We have the canonical
morphism
$$\can:\ \Dd= A\ot_B A\to \Cc=\bigoplus_{\sigma\in G} Av_\sigma,$$
given by
$$\can(a\ot b)= \sum_{\sigma\in G} av_\sigma b=
\sum_{\sigma\in G} a\alpha_\sigma(be_{\sigma^{-1}})v_\sigma$$
We can also compute that
$${}^*\can:\ {}^*\Cc=\bigoplus_{\sigma\in G}u_\sigma A\to {}^*\Dd\cong
{}_B\End(A)^{\rm op}$$
is given by
$${}^*\can(u_\tau b_\tau)(a)=\alpha_\tau(ae_{\tau^{-1}})b_\tau.$$
Let us now compute the module $Q\subset {}^*\Cc$ introduced in \seref{1.1}
Recall that $q\in Q$ if and only if
\begin{equation}\eqlabel{4.8.1}
c_{(1)}q(c_{(2)})=q(c)x,
\end{equation}
for all $c\in \Cc$.

\begin{proposition}\prlabel{4.8}
$Q=\{\sum_{\sigma\in G} u_\sigma\alpha_\sigma(ae_{\sigma^{-1}})~|~a\in A\}$.
\end{proposition}

\begin{proof}
Take $q=\sum_{\sigma\in G} u_\sigma a_\sigma\in Q$, with $a_\sigma\in Ae_\sigma$,
and put $c=v_\tau$ in
\equref{4.8.1}. Recall that $\Delta_\Cc(v_\sigma)=\sum_{\rho\in G}
v_\rho\ot_A v_{\rho^{-1}\tau}$. Then we calculate that
$$
c_{(1)}q(c_{(2)})=\sum_{\rho,\sigma\in G} v_\rho\delta_{\sigma,\rho^{-1}\tau}
a_{\sigma}=
\sum_{\rho\in G} v_\rho a_{\rho^{-1}\tau}=
\sum_{\rho \in G} \alpha_\rho(a_{\rho^{-1}\tau}e_{\rho^{-1}})v_\rho,
$$
$$q(c)=(\sum_{\sigma\in G} u_\sigma a_\sigma)(v_\tau)=a_\tau$$
and
$$q(c)x=\sum_{\rho \in G} a_\tau v_\rho$$
Hence it follows that
\begin{equation}\eqlabel{4.8.2}
a_\tau e_\rho= \alpha_\rho(a_{\rho^{-1}\tau}e_{\rho^{-1}}),
\end{equation}
for all $\tau,\rho\in G$. Taking $\tau=\rho$, we find that
\begin{equation}\eqlabel{4.8.3}
a_\tau= a_\tau e_\tau=\alpha_\tau(a_1e_{\tau^{-1}}),
\end{equation}
and we find that
\begin{equation}\eqlabel{4.8.4}
q=\sum_{\sigma\in G} u_\sigma a_\sigma=
\sum_{\sigma\in G} u_\sigma\alpha_\sigma(a_1e_{\sigma^{-1}})
\end{equation}
is of the desired form. Conversely, take $q$ of the form \equref{4.8.4}.
Then \equref{4.8.3} holds. Using \equref{4.2.1}, we compute
$$
\alpha_\rho(a_{\rho^{-1}\tau}e_{\rho^{-1}})=
\alpha_{\rho}(\alpha_{\rho^{-1}\tau}(a_1e_{\tau^{-1}\rho})e_{\rho^{-1}})=
 \alpha_\tau(a_1e_{\tau^{-1}})e_\rho= a_\tau e_\rho,$$
and \equref{4.8.2} follows, which means that \equref{4.8.1} holds for
$c=v_\tau$. Using the left $A$-linearity of $q$ and $\Delta_\Cc$, it follows
that \equref{4.8.1} holds for arbitrary $c\in \Cc$.
\end{proof}

It follows from \prref{4.8} that we have an isomorphism of abelian groups
$$A\to Q,~~a\mapsto \sum_{\sigma\in G} u_\sigma\alpha_\sigma(ae_{\sigma^{-1}}).$$
This can also be seen using \prref{4.7} and \cite[Theorem 2.7]{CVW}.
The $({}^*\Cc,T)$-bimodule structure on $Q$ (see \prref{2.2}) can be transported
to $A$. The right $T$-action on $A$ is then given by right multiplication,
and the left ${}^*\Cc$-action is the following:
$$(u_\tau a_\tau)\cdot a=\alpha_{\tau^{-1}}(a_\tau ae_\tau).$$
Recall also from \prref{2.2} that $A\in {}_T\Mm_{{}^*\Cc}$. The left $T$-action
is given by left multiplication. The right ${}^*\Cc$-action is the following:
\begin{eqnarray*}
&&\hspace*{-2cm}
a\cdot(u_\tau a_\tau)=(au_\tau a_\tau)(x)\\
&=& \sum_{\sigma\in G} u_\tau \alpha_\tau(ae_{\tau^{-1}})a_\tau)(v_\sigma)
=\alpha_\tau(a e_{\tau^{-1}})b_\tau.
\end{eqnarray*}

We have seen in \prref{2.2} that we have a Morita context
$(T,{}^*\Cc,A,Q,\tau,\mu)$. Using the isomorphism between $A$ and $Q$,
we find a Morita context $(T,{}^*\Cc,A,A,\tau,\mu)$. Let us compute the
connecting maps $\tau:\ A\ot_{{}^*\Cc}A\to T$ and
$\mu:\  A\ot_{T}A\to {}^*\Cc$, using (\ref{eq:2.2.1}-\ref{eq:2.2.2}).
\begin{eqnarray*}
&&\hspace*{-2cm}\tau(b\ot a)=
(\sum_{\sigma\in G}u_\sigma \alpha_\sigma(ae_{\sigma^{-1}}))
(\sum_{\tau\in G} (v_\tau b))\\
&=&\sum_{\sigma,\tau\in G} u_{\sigma}(\alpha_\tau(be_{\tau^{-1}})v_\tau)
\alpha_\sigma(ae_{\sigma^{-1}})\\
&=& \sum_{\sigma\in G} \alpha_\sigma(be_{\sigma^{-1}})\alpha_\sigma(ae_{\sigma^{-1}})\\
&=& \sum_{\sigma\in G} \alpha_\sigma(bae_{\sigma^{-1}});\\
\mu(a\ot b)&=& \sum_{\sigma\in G}u_\sigma \alpha_\sigma(ae_{\sigma^{-1}})b.
\end{eqnarray*}
We summarize our results as follows.

\begin{proposition}\prlabel{4.9}
We have a Morita context $(T,{}^*\Cc,A,Q,\tau,\mu)$. The connecting maps are
given by the formulas
\begin{eqnarray}
\tau(b\ot a)&=& \sum_{\sigma\in G} \alpha_\sigma(bae_{\sigma^{-1}});\eqlabel{4.9.1}\\
\mu(a\ot b)&=& \sum_{\sigma\in G}u_\sigma \alpha_\sigma(ae_{\sigma^{-1}})b\eqlabel{4.9.2}.
\end{eqnarray}
\end{proposition}

\begin{proposition}\prlabel{4.10}
The map $\tau$ in the Morita context from \prref{4.9} is surjective if and only if
there exists $a\in A$ such that
$$\sum_{\sigma\in G} \alpha_\sigma(ae_{\sigma^{-1}})=1.$$
\end{proposition}

\begin{proof}
According to \cite[Theorem 3.3]{CVW}, $\tau$ is surjective if and only if there exists $q\in Q$
such that $q(x)=1$. Let $a\in A$ correspond to $q\in Q$. Then we compute that
$$q(x)=(\sum_{\sigma\in G} u_\sigma\alpha_\sigma(ae_{\sigma^{-1}}))
(\sum_{\tau\in G} v_\tau)= \sum_{\sigma\in G} \alpha_\sigma(ae_{\sigma^{-1}}),$$
and the result follows.
\end{proof}

From \thref{2.4}, we obtain:

\begin{theorem}\thlabel{4.11}
Let $G$ be a finite group, and $(e_\sigma,\alpha_\sigma)_{\sigma\in G}$ an
idempotent partial action of $G$ on $A$. Let $i:\ B\to T=A^{{\rm co}\Cc}$ a ring
morphism, and consider $\can:\ A\ot_BA\to \Cc$. Then the following assertions
are equivalent. 
\begin{enumerate}
\item \begin{itemize}
\item ${\rm can}$ is an isomorphism;
\item $A$ is faithfully flat as a left $B$-module.
\end{itemize}
\item \begin{itemize}
\item ${}^*{\rm can}$ is an isomorphism;
\item $A$ is a left $B$-progenerator.
\end{itemize}
\item \begin{itemize}
\item $B=T$;
\item the Morita context $(B,{}^*\Cc,A,A,\tau,\mu)$ is strict.
\end{itemize}
\item \begin{itemize}
\item $B=T$;
\item $(F,G)$ is an equivalence of categories.
\end{itemize}
\end{enumerate}
\end{theorem}

If we take $A$ and $B$ commutative, then \thref{4.11} implies part of
\cite[Theorem 3.1]{DFP}, namely the equivalence of the conditions (i),
(ii) and (iii).

\end{document}